\documentclass[a4paper,12pt,twoside]{article}
\usepackage{amsmath,amssymb}
\usepackage[amsmath,thmmarks]{ntheorem}
\usepackage{graphicx}
\usepackage{new}
\theoremseparator{.}
\theorembodyfont{\upshape}
\newtheorem{defn}{Definition}[section]
\newtheorem{rem}[defn]{Remark}
\theorembodyfont{\slshape}
\newtheorem{thm}[defn]{Theorem}
\newtheorem{prop}[defn]{Proposition}
\newtheorem{lem}[defn]{Lemma}

\theoremstyle{nonumberplain}
\theorembodyfont{\upshape}
\theoremsymbol{\rule{.4em}{.8em}}
\theoremheaderfont{\itshape}
\newtheorem{prf}{Proof}
\def\whatistobeproved{\relax}
\makeatletter
\newtheorem{@prfof}{Proof of \whatistobeproved}
\def\prfof#1{\def\whatistobeproved{#1}\@prfof}
\def\endprfof{\end@prfof}
\makeatother
\newcommand{\word}[1]{\textit{#1}}

\DeclareMathOperator{\diam}{diam}
\newcommand{\B}{\mathcal{B}}
\newcommand{\Bdd}{\boldsymbol{b}}
\newcommand{\Q}{\mathbb{Q}}
\newcommand{\R}{\mathbb{R}}
\newcommand{\add}{\boldsymbol{\Sigma}}

\newcommand{\amb}{\boldsymbol{\Delta}}
\title{Continuity Points of Typical Bounded Functions}
\author{Shingo Saito}
\begin{document}
\begin{affiliation}
 Faculty of Mathematics,
 Kyushu University,
 6--10--1, Hakozaki, Higashi-ku, Fukuoka, 812--8581, Japan.

 Email address: \texttt{ssaito@math.kyushu-u.ac.jp}


\end{affiliation}

\maketitle

\begin{abstract}
 Kostyrko and \v{S}al\'{a}t showed that
 if a linear space of bounded functions has an element
 that is discontinuous almost everywhere,
 then a typical element in the space is discontinuous almost everywhere.
 We give a topological analogue of this theorem and provide some examples.
\end{abstract}

\section{Introduction}
The study of typical continuous functions
has been one of the most popular topics in classical real analysis.
Several people have also investigated typical behaviour in other families of functions,
such as those of bounded functions of Baire class~$1$,
of bounded Darboux functions of Baire class~$1$, and
of bounded derivatives,
where the topology is given by the supremum norm.

Kostyrko and \v{S}al\'{a}t \cite{MR0795616} proved
a theorem, applicable to general families of bounded functions,
on the continuity points of typical functions.
They showed that if a linear space of bounded functions has an element
that is discontinuous almost everywhere, then a typical element in the space
is discontinuous almost everywhere
(what they actually proved is slightly stronger than this; see below for the precise statement).
Our aim in this paper is to give a topological analogue of this theorem,
by showing that if a linear space of bounded functions has an element
that is discontinuous everywhere in a residual set, then a typical element in the space
is discontinuous everywhere in a residual set
(again, we in fact show a slightly stronger result).

Let us fix notation and give the precise statements
of the theorem of Kostyrko and \v{S}al\'{a}t and of our main theorem.
By a \word{function} we shall always mean a real-valued function defined on the unit interval $[0,1]$.
Let $\Bdd$ denote the Banach space of all bounded functions, equipped with the supremum norm.
For a function $f$, we write $C(f)$ and $D(f)$
for the sets of all continuity and discontinuity points of $f$ respectively.
The Lebesgue measure on $[0,1]$ will be denoted by $\mu$.

\begin{thm}[{{\cite[Theorem]{MR0795616}}}]
 If $X$ is a linear subspace of $\Bdd$ such that
 \[
  \inf_{f\in X}\mu\bigl(C(f)\bigr)=0,
 \]
 then the set
 \[
  \bigl\{f\in X\bigm|\mu\bigl(C(f)\bigr)=0\bigr\}
 \]
 is a residual $G_{\delta}$ subset of $X$.
\end{thm}

\begin{rem}
 The hypothesis is clearly fulfilled if there exists $f\in X$ for which $\mu\bigl(C(f)\bigr)=0$.
 Note also that the residual $G_{\delta}$ subset in the conclusion may be empty
 because $X$ is not assumed to be a closed subspace of $\Bdd$.
\end{rem}

Our topological analogue is the following:

\begin{thm}[Main Theorem]\label{thm:cont_mgr}
 Let $X$ be a linear subspace of $\Bdd$ such that
 for each nonempty open subset $U$ of $[0,1]$, there exists $f\in X$ such that $C(f)$ is not residual in $U$.
 Then the set
 \[
  \{f\in X\mid\text{$C(f)$ is nowhere dense}\}
 \]
 is a residual $G_{\delta}$ subset of $X$.
\end{thm}

\section{Proof of the main theorem}
\begin{defn}
 Let $f$ be a function.
 For each $x\in[0,1]$, we define the \word{oscillation} $\omega(f,x)$ of $f$ at $x$ by
 \[
  \omega(f,x)=\inf\{\diam f(U)\mid\text{$U$ is an open neighbourhood of $x$}\}\in[0,\infty],
 \]
 where $\diam f(U)$ denotes the diameter of $f(U)$.
 For each $t\in\R$, we write
 \[
  D_t(f)=\{x\in[0,1]\mid\omega(f,x)\ge t\}.
 \]
\end{defn}

We denote by $\Q_+$ the set of all positive rational numbers.

\begin{lem}\label{lem:D(f)_F_sigma}
 Let $f$ be a function.
 Then $D_t(f)$ is closed for every $t\in\R$, and we have $\bigcup_{t>0}D_t(f)=\bigcup_{t\in\Q_+}D_t(f)=D(f)$.
 It follows that $D(f)$ is $F_{\sigma}$ and $C(f)$ is $G_{\delta}$.
\end{lem}

\begin{prf}
 Obvious (see \cite[Theorem~7.1]{MR0584443} for details).
\end{prf}

\begin{lem}\label{lem:D_estimate}
 If $f$ and $g$ are functions and $t\in\R$, then $D_t(f)\subset D_{t-2\lVert f-g\rVert}(g)$.
\end{lem}

\begin{prf}
 Let $x\in D_t(f)$.
 Given $\varepsilon>0$ and a neighbourhood $U$ of $x$, we have
 \[
  \diam f(U)\ge\omega(f,x)\ge t>t-\varepsilon,
 \]
 which implies that there exist $x_1,x_2\in U$ with $\lvert f(x_1)-f(x_2)\rvert>t-\varepsilon$.
 It follows that
 \begin{align*}
  \diam g(U)&\ge\lvert g(x_1)-g(x_2)\rvert
  \ge\lvert f(x_1)-f(x_2)\rvert-2\lVert f-g\rVert\\
  &>t-2\lVert f-g\rVert-\varepsilon.
 \end{align*}
 Taking the infimum over all neighbourhoods $U$, we obtain
 $\omega(g,x)\ge t-2\lVert f-g\rVert-\varepsilon$.
 Since $\varepsilon$ was arbitrary, we see that $\omega(g,x)\ge t-2\lVert f-g\rVert$,
 i.e.~$x\in D_{t-2\lVert f-g\rVert}(g)$.
\end{prf}

Let $X$ be a linear subspace of $\Bdd$ satisfying the assumption of Theorem~\ref{thm:cont_mgr}.

In what follows, by a \word{subinterval} we shall always mean
a nondegenerate closed subinterval of $[0,1]$ with rational endpoints.
For each subinterval $I$,
we write
\begin{align*}
 X_I&=\{f\in X\mid\text{$C(f)$ is not residual in $I$}\}\\
 &=\{f\in X\mid\text{$D(f)$ is nonmeagre in $I$}\}.
\end{align*}

\begin{lem}\label{lem:X_I_open}
 The set $X_I$ is open in $X$ for every subinterval $I$.
\end{lem}

\begin{prf}
 Let $f\in X_I$ be given.
 Since $D(f)=\bigcup_{t\in\Q_+}D_t(f)$ is nonmeagre in $I$,
 there exist $t\in\Q_+$ and a subinterval $J\subset I$ such that $D_t(f)\supset J$.
 If $g\in X$ satisfies $\lVert g-f\rVert<t/2$, then
 Lemma~\ref{lem:D_estimate} shows that
 \[
  J\subset D_t(f)\subset D_{t-2\lVert g-f\rVert}(g)\subset D(g),
 \]
 and so $g\in X_I$.
\end{prf}

\begin{lem}\label{lem:X_I_dense}
 The set $X_I$ is dense in $X$ for every subinterval $I$.
\end{lem}

\begin{prf}
 Given $f\in X$ and $\varepsilon>0$,
 we need to find $h\in X_I$ with $\lVert f-h\rVert<\varepsilon$.
 We may assume that $f\notin X_I$, i.e.~$C(f)$ is residual in $I$.
 By the assumption on $X$, there exists $g\in X$ for which
 $C(g)$ is not residual in $I$.
 Choose $c>0$ so small that $c\lVert g\rVert<\varepsilon$,
 and set $h=f+cg$.
 Observe that $h\in X$ because $X$ is a linear space.
 Now $C(h)$ cannot be residual in $I$ since $C(f)\cap C(h)\subset C(g)$.
\end{prf}

\begin{prfof}{Theorem~\ref{thm:cont_mgr}}
 By Lemmas \ref{lem:X_I_open} and \ref{lem:X_I_dense},
 it suffices to show that $C(f)$ is nowhere dense
 if and only if it is not residual in any subinterval.
 If $C(f)$ is nowhere dense, then in any subinterval it is not dense and so not residual.
 If $C(f)$ is not nowhere dense, then in some subinterval it is dense $G_{\delta}$ by Lemma~\ref{lem:D(f)_F_sigma}
 and so residual.
\end{prfof}

\section{Examples}
If a linear subspace $X$ of $\Bdd$ consists entirely of functions of Baire class~$1$,
then $C(f)$ is residual for every $f\in X$ (see \cite[Theorem~7.3]{MR0584443}), and so
$X$ does not satisfy the assumption of our main theorem.
Let us consider the family $X=\Bdd\B_2$ of all bounded functions of Baire class~$2$.
Since the characteristic function of $\Q$ is a nowhere continuous function of Baire class~$2$,
the theorem shows that $C(f)$ is nowhere dense for a typical $f\in\Bdd\B_2$.
In fact, it turns out that a much stronger result holds:

\begin{prop}\label{prop:typical_Baire_alpha_nowhere_cont}
 If $\alpha$ is an ordinal with $2\le\alpha<\omega_1$
 and $\Bdd\B_{\alpha}$ denotes the family of all bounded functions of Baire class~$\alpha$,
 then $C(f)=\emptyset$ for a typical $f\in\Bdd\B_{\alpha}$.
 In other words, a typical $f\in\Bdd\B_{\alpha}$ is nowhere continuous.
\end{prop}

\begin{rem}
 The space $\Bdd\B_{\alpha}$ is closed in $\Bdd$ and therefore complete.
\end{rem}

\begin{prfof}{Proposition~\ref{prop:typical_Baire_alpha_nowhere_cont}}
 Note first that the set
 \[
  \{f\in\Bdd\B_{\alpha}\mid\text{$D_t(f)=[0,1]$ for some $t>0$}\}
 \]
 is open because of Lemma~\ref{lem:D_estimate}.
 Consequently, it suffices to prove that the set is dense as well.

 Let $f\in\Bdd\B_{\alpha}$ and $\varepsilon>0$ be given.
 Take $a,b\in\R$ with $f([0,1])\subset(a,b)$,
 and choose an integer $n\ge2$ so large that $2(b-a)/n<\varepsilon$.
 Put $h=(b-a)/n$.

 For each $k=0,\dots,n-2$, set
 \[
  A_k=\{x\in[0,1]\mid a+kh<f(x)<a+(k+2)h\}.
 \]
 Then $A_k\in\add^0_{\alpha+1}$ for every $k$ and $\bigcup_{k=0}^{n-2}A_k=[0,1]$.
 By the reduction property of $\add^0_{\alpha+1}$ (see \cite[Theorem~22.16]{MR1321597}),
 we may find disjoint sets $B_k\in\add^0_{\alpha+1}$ such that
 $B_k\subset A_k$ for every $k$ and $\bigcup_{k=0}^{n-2}B_k=[0,1]$.

 Define a function $g$ by
 \[
  g(x)=
  \begin{cases}
   a+kh&\text{if $x\in B_k\cap\Q$};\\
   a+(k+1/2)h&\text{if $x\in B_k\setminus\Q$}.
  \end{cases}
 \]
 Note that $g\in\Bdd\B_{\alpha}$
 because $\Q\in\add^0_2\subset\amb^0_{\alpha+1}$ due to $\alpha\ge2$.
 Observe that $\lVert f-g\rVert<2h<\varepsilon$.

 Given any $x\in[0,1]$ and a neighbourhood $U$ of $x$,
 we shall find $y\in U$ with $\lvert g(y)-g(x)\rvert\ge h/2$;
 this implies that $D_{h/2}(g)=[0,1]$ and finishes the proof.
 Take the unique $k$ with $x\in B_k$.
 If there exists $y\in U\setminus B_k$, then $\lvert g(y)-g(x)\rvert\ge h/2$.
 Otherwise, choosing $y\in U$ so that exactly one of $x$ and $y$ is rational,
 we obtain $\lvert g(y)-g(x)\rvert=h/2$.
\end{prfof}

The same method also works for some other families:

\begin{prop}
 \begin{enumerate}
  \item A typical $f\in\Bdd$ is nowhere continuous.
  \item A typical Lebesgue measurable $f\in\Bdd$ is nowhere continuous.
  \item A typical $f\in\Bdd$ with the Baire property is nowhere continuous.
 \end{enumerate}
\end{prop}

\begin{prf}
 Mimic the proof of Proposition~\ref{prop:typical_Baire_alpha_nowhere_cont}.
\end{prf}

\begin{rem}
 Statement (3) is clearly stronger than Corollary~IV in \cite{MR0795616},
 which asserts that in the space of all bounded functions with the Baire property,
 the functions discontinuous almost everywhere form a residual $G_{\delta}$ subset.
\end{rem}

Now we give an example of $X$
for which the assumption of our main theorem holds
but none of whose members is nowhere continuous:

\begin{prop}
 If $A$ is a nowhere dense $G_{\delta}$ subset of $[0,1]$, then
 the family
 \[
  X=\{f\in\Bdd\mid C(f)\supset A\}
 \]
 satisfies the hypothesis of Theorem~\ref{thm:cont_mgr}.
\end{prop}

\begin{prf}
 Remember that every $G_{\delta}$ set is the set of continuity points of some $f\in\Bdd$
 (see \cite[Theorem~7.2]{MR0584443}).
\end{prf}

\begin{rem}
 The space $X$ is a closed subspace because if $f_n\in X$ and $f_n\to f$ uniformly,
 then $C(f)\supset\bigcap_{n=1}^{\infty}C(f_n)\supset A$.
\end{rem}

\end{document}